\documentclass[11pt]{article}

\usepackage{amsfonts,epsf,amsmath,tikz}

\newtheorem{thm}{Theorem}[section]
\newtheorem{prp}[thm]{Proposition}

\newtheorem{conj}[thm]{Conjecture}

\newtheorem{lemma}[thm]{Lemma}

\newcommand{\qed}{\hfill ~$\square$\bigskip}
\newcommand{\proof}{\noindent{\bf Proof.} }

\newcommand{\cp}{\,\square\,}

\newcommand{\D}{{Dominator }}
\newcommand{\St}{{Staller }}
\begin{document}

\title{Domination game played on trees and spanning subgraphs}
\author{
Bo\v{s}tjan Bre\v{s}ar\thanks{Supported by the
Ministry of Science of Slovenia
under the grants P1-0297 and J1-2043. The author is also with the
Institute of Mathematics, Physics and Mechanics,
Jadranska 19, 1000 Ljubljana.} \\
Faculty of Natural Sciences and Mathematics \\
University of Maribor, Slovenia\\
bostjan.bresar@uni-mb.si
\and
Sandi Klav\v zar$^{\ast}$ \\
Faculty of Mathematics and Physics\\
University of Ljubljana, Slovenia \\
and\\
Faculty of Natural Sciences and Mathematics \\
University of Maribor, Slovenia\\
sandi.klavzar@fmf.uni-lj.si
\and
Douglas F. Rall\thanks{
Supported by the Wylie Enrichment Fund of Furman University
and by a grant from the Simons Foundation (Grant Number 209654 to Douglas F. Rall).
Part of the research done during a sabbatical visit at the University of Ljubljana
and University of Maribor.}
\\
Herman N. Hipp Professor of Mathematics \\
Department of Mathematics, Furman University\\
Greenville, SC, USA\\
doug.rall@furman.edu
}

\date{\today}

\maketitle

\begin{abstract}
The domination game, played on a graph $G$, was introduced in~\cite{brklra-2010}.
Vertices are chosen, one at a time, by two players \D and Staller. Each chosen
vertex must enlarge the set of vertices of $G$ dominated to that point in the game.
Both players use an optimal strategy---Dominator plays so as to end the game as quickly
as possible, and Staller plays in such a way that the game lasts as many steps as possible.
The game domination number $\gamma_g(G)$ is the number of
vertices chosen when \D starts the game and the Staller-start game
domination number $\gamma_g'(G)$ is the result when \St starts the game.

In this paper these two games are studied when played on trees and spanning subgraphs.
A lower bound for the game domination number
of a tree in terms of the order and maximum degree is proved and shown
to be asymptotically tight. It is shown that for every $k$, there is
a tree $T$ with $(\gamma_g(T),\gamma_g'(T)) = (k,k+1)$ and conjectured
that there is none with $(\gamma_g(T),\gamma_g'(T)) = (k,k-1)$.
A relation between the game domination number of a graph and its spanning
subgraphs is considered. It is proved that there exist 3-connected graphs
$G$ having a 2-connected spanning subgraph $H$ such that the game
domination number of $H$ is arbitrarily smaller than that of $G$.
Similarly, for any integer $\ell \geq 1$, there exists a graph $G$
and a spanning tree $T$ such that $\gamma_g(G)-\gamma_g(T)\ge \ell$.
On the other hand, there exist graphs $G$ such that the game domination number
of any spanning tree of $G$ is arbitrarily larger than that of $G$.
\end{abstract}



\section{Introduction}

The domination game played on a graph $G$ consists of two players,
\D and Staller, who alternate taking turns choosing a vertex from $G$
such that whenever a vertex is chosen by either player, at least one
additional vertex is dominated. \D wishes to dominate the graph as fast
as possible, and \St wishes to delay the process as much as possible.
The {\em game domination number}, denoted $\gamma_g(G)$, is the number of
vertices chosen when \D starts the game provided that both
players play optimally. Similarly,  the {\em Staller-start game
domination number}, written as $\gamma_g'(G)$, is the result of the game when \St starts the game.
The Dominator-start game and the Staller-start game will be briefly
called {\em Game 1} and {\em Game 2}, respectively.

This game was first studied in 2010 (\cite{brklra-2010}) but was brought to the authors'
attention back in 2003 by Henning~\cite{he-03}. Among other results, the
authors of~\cite{brklra-2010} proved a lower bound for the game domination
number of the Cartesian product of graphs and established a connection
with Vizing's conjecture; for the latter see~\cite{brdo-2012}.
The Cartesian product was further investigated in~\cite{bill-2011}, where the behavior
of $\lim_{\ell \to \infty} \gamma_g(K_m \cp P_\ell)/\ell$ was studied in detail.

In the rest of this section we give some notation and definitions, and we recall results
needed later. In Section~\ref{sec:lower} we prove a general lower bound for the game domination
number of a tree. In Section~\ref{sec:pairs} we consider which pairs of integers $(r,s)$
can be realized as $(\gamma_g(T),\gamma_g'(T))$, where $T$ is a tree.
It is shown that this is the case for all pairs but those
of the form $(k,k-1)$. This enlarges the set of pairs known to be realizable by {\em connected} graphs.
We conjecture that the pairs $(k,k-1)$ cannot be realized by trees. In the final section we study
relations between the game domination number of a graph and its spanning
subgraphs. We construct graphs $G$ having a spanning tree $T$
with $\gamma_g(G)-\gamma_g(T)$ arbitrarily large and build 3-connected graphs having a 2-connected
spanning subgraph exhibiting the same phenomenon. This is rather surprising
because the domination number of a spanning tree (or a spanning subgraph)
can never be smaller than the domination number of its supergraph.
We also present graphs $G_{2r}$ ($r\ge 1$) such that $\gamma_g(T) - \gamma_g(G_{2r}) \ge r-1$
holds for any spanning tree $T$ of $G_{2r}$. This is again different from the usual domination
because it is known (see~\cite[Exercise 10.14]{imklra-2008}) that every graph contains a spanning
tree with the same domination number.

Throughout the paper we will use the convention that $d_1, d_2, \ldots$
denotes the list of vertices chosen by \D and $s_1,s_2, \ldots$ the
list chosen by Staller.
We say that a pair $(r,s)$ of integers is {\em realizable} if there exists
a graph $G$ such that $\gamma_g(G) = r$ and $\gamma_g'(G) = s$.
Following~\cite{bill-2011}, we make the following definitions.
A {\em partially dominated graph} is a graph in which some vertices have already
been dominated in some turns of the game already played. A vertex $u$ of a
partially dominated graph $G$ is {\em saturated} if each vertex in $N[u]$ is
dominated. The {\em residual graph} of $G$ is the graph obtained from $G$
by removing all saturated vertices and all edges joining dominated vertices.
If $G$ is a partially dominated graph, then $\gamma_g(G)$ and $\gamma_g'(G)$ denote the
optimal number of moves remaining in Game 1 and Game 2, respectively
(it is assumed here that Game 1, respectively Game 2, refers to Dominator,
respectively Staller, being the first to play in the partially dominated graph $G$).

The game domination number of a graph $G$ can be bounded in terms of the domination number
$\gamma(G)$ of $G$:

\begin{thm}
\label{thm:connection}
{\rm (\cite{brklra-2010})}
For any graph $G$, $\gamma(G) \le \gamma_g(G) \le 2\gamma(G)-1$.
\end{thm}

It was demonstrated in~\cite{brklra-2010} that, in general, Theorem~\ref{thm:connection}
cannot be improved. More precisely, for any positive integer $k$ and any integer $r$ such
that $0\le r\le k-1$, there exists a graph $G$ with $\gamma(G) = k$ and $\gamma_g(G) = k+r$.

The game domination number and the Staller-start game domination
number never differ by more than 1 as the next result asserts.

\begin{thm}
{\rm (\cite{brklra-2010,bill-2011})}
\label{thm:we-and-bill}
If $G$ is any graph, then $|\gamma_g(G) - \gamma_g'(G)| \le 1$.
\end{thm}

By Theorem~\ref{thm:we-and-bill} only pairs of the form $(r,r)$,
$(r,r+1)$, and $(r,r-1)$ are realizable. See~\cite{ko-2014} for
a study of realizable pairs.

The following lemma, due to Kinnersley, West, and Zamani~\cite{bill-2011}
in particular implies $\gamma_g'(G) \le \gamma_g(G) + 1$, which is one
half of Theorem~\ref{thm:we-and-bill}. The other half was earlier
proved in~\cite{brklra-2010}.

\begin{lemma}
{\rm (Continuation Principle)}
\label{lem:continuation}
Let $G$ be a graph and let $A$ and $B$ be subsets of $V(G)$. Let $G_A$ and $G_B$ be partially dominated
graphs in which the sets $A$ and $B$ have already been dominated, respectively.
If $B\subseteq A$, then $\gamma_g(G_A) \le \gamma_g(G_B)$ and
$\gamma_g'(G_A) \le \gamma_g'(G_B)$.
\end{lemma}

We wish to point out that the Continuation Principle is a very useful tool for
proving results about game domination number. In particular, suppose that
at some stage of the game a vertex $x$ is an optimal move for Dominator.
If for some vertex $y$ such that the undominated part of $N[x]$ is
contained in $N[y]$, then $y$ is also an optimal selection for Dominator.
We can thus assume (if desired) that he plays $y$.

\section{A lower bound for trees}
\label{sec:lower}

In this section we give a lower bound on the game domination number of
trees and show that it is asymptotically sharp. Before we can state the main
result, we need the following:

\begin{lemma}
\label{lem:key-for-trees}
In a partially dominated tree $F$, \St can
make a move that dominates at most two new vertices.
\end{lemma}

\proof
Let $A$ be the set of saturated vertices of $F$ and let
$B$ be the set of vertices of $F$ that are dominated but not saturated.
Let $C = V(F) - (A\cup B)$. Let $F'$ be the subforest of $F$ induced
by $B\cup C$ but with edges induced by $B$ removed (that is, $F'$ is the
residual graph). We may assume that $C\ne \emptyset$.
Now $F'$ has a leaf $x$. If \St plays $x$, then she dominates
at most two vertices in $C$. If $x\in B$, then
\St dominates exactly one vertex in $C$.
\qed

The move guaranteed by Lemma~\ref{lem:key-for-trees} may not be
an optimal move for Staller. For instance, the optimal first move of
\St when playing on $P_5$ is the middle vertex of $P_5$, thus
dominating three new vertices. Also, we will see later that an
optimal first move for \St when playing Game 2 on the tree $T_r$ from
Figure~\ref{fig:3realizations} is $w$, thus dominating $r+1$ new
vertices.

\begin{thm}
\label{thm:lower-for-trees}
If $T$ is a tree on $n$ vertices, then
$$\gamma_g(T) \ge \left\lceil \frac{2n}{\Delta(T)+3} \right\rceil -1\,.$$
\end{thm}

\proof
By Lemma~\ref{lem:key-for-trees}, \St can move in such a way that at most
two new vertices are dominated on each of her moves. Let us suppose that
\D plays optimally when \St plays to dominate at most two new vertices on
each move. Let $d_1, s_1, d_2, s_2, \ldots, d_t, s_t$ be the resulting game,
where we assume that $s_t$ is the empty move if $T$ is dominated
after the move $d_t$. Let $f(d_i)$ (resp. $f(s_i)$) denote the number of
newly dominated vertices when \D plays $d_i$ (resp. when \St plays $s_i$).
If the game ends on a move by Staller, then
$$n = \sum_{i=1}^t \left( f(d_i) + f(s_i) \right)
\le \sum_{i=1}^t \left( ({\rm deg}(d_i)+1) +2 \right)
= \sum_{i=1}^t {\rm deg}(d_i)+3t\,.$$
Since this strategy may not be an optimal one for Staller, it follows that
$\gamma_g(T)\ge 2t$. Similar arguments give $\gamma_g(T) \ge 2t-1$ if the game
ends on Dominator's move.

If \St ends the game, then $n \le t(\Delta(T) + 3) \le \frac{1}{2}\gamma_g(T)(\Delta(T) + 3)$,
and hence $\gamma_g(T) \ge \left\lceil \frac{2n}{\Delta(T)+3} \right\rceil$
since $\gamma_g(T)$ is integral. If the game ends on Dominator's move,
then $\gamma_g(T) \ge 2t-1$, and hence
$$n \le t(\Delta(T)+3)-2 \le t(\Delta(T) + 3) \le \frac{\gamma_g(T)+1}{2}(\Delta(T) + 3)\,.$$
This is equivalent to $2n\le (\gamma_g(T)+1)(\Delta(T)+3)$, which in turn implies that
$$\gamma_g(T) \ge \left\lceil \frac{2n}{\Delta(T)+3} -1 \right\rceil
= \left\lceil \frac{2n}{\Delta(T)+3} \right\rceil - 1\,,$$
as claimed.
\qed

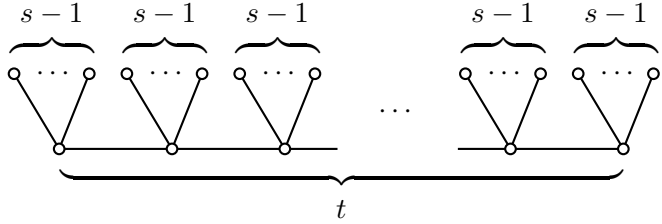
\begin{figure}[ht!]
\begin{center}
\begin{tikzpicture}[scale=1.0,style=thick]
\def\vr{2pt} 
\draw (0.0,0.0) -- (1.5,0.0); \draw (1.5,0.0) -- (3.0,0.0);
\draw (3.0,0.0) -- (3.7,0.0); \draw (5.3,0.0) -- (6.0,0.0);
\draw (6.0,0.0) -- (7.5,0.0);
\draw (0.0,0.0) -- (-0.6,1.0);
\draw (0.0,0.0) -- (0.4,1.0);
\draw (1.5,0.0) -- (0.9,1.0);
\draw (1.5,0.0) -- (1.9,1.0);
\draw (3.0,0.0) -- (2.4,1.0);
\draw (3.0,0.0) -- (3.4,1.0);
\draw (6.0,0.0) -- (5.4,1.0);
\draw (6.0,0.0) -- (6.4,1.0);
\draw (7.5,0.0) -- (6.9,1.0);
\draw (7.5,0.0) -- (7.9,1.0);
\draw (0.0,0.0)  [fill=white] circle (\vr);
\draw (1.5,0.0)  [fill=white] circle (\vr);
\draw (3.0,0.0)  [fill=white] circle (\vr);
\draw (6.0,0.0)  [fill=white] circle (\vr);
\draw (7.5,0.0)  [fill=white] circle (\vr);
\draw (-0.6,1.0)  [fill=white] circle (\vr);
\draw (0.4,1.0)  [fill=white] circle (\vr);
\draw (0.9,1.0)  [fill=white] circle (\vr);
\draw (1.9,1.0)  [fill=white] circle (\vr);
\draw (2.4,1.0)  [fill=white] circle (\vr);
\draw (3.4,1.0)  [fill=white] circle (\vr);
\draw (5.4,1.0)  [fill=white] circle (\vr);
\draw (6.4,1.0)  [fill=white] circle (\vr);
\draw (6.9,1.0)  [fill=white] circle (\vr);
\draw (7.9,1.0)  [fill=white] circle (\vr);
%
\draw (3.75,-0.3) node {$\underbrace{\phantom{xxxxxxxxxxxxxxxxxxxxxxxxxxxxxxxxxx}}$};
\draw (3.75,-0.8) node {$t$};
\draw (-0.05,1.0) node {$\cdots$};
\draw (-0.1,1.3) node {$\overbrace{\phantom{xxxxx}}$};
\draw (-0.1,1.8) node {$s-1$};
\draw (1.45,1.0) node {$\cdots$};
\draw (1.4,1.3) node {$\overbrace{\phantom{xxxxx}}$};
\draw (1.4,1.8) node {$s-1$};
\draw (2.95,1.0) node {$\cdots$};
\draw (2.9,1.3) node {$\overbrace{\phantom{xxxxx}}$};
\draw (2.9,1.8) node {$s-1$};
\draw (5.95,1.0) node {$\cdots$};
\draw (5.9,1.3) node {$\overbrace{\phantom{xxxxx}}$};
\draw (5.9,1.8) node {$s-1$};
\draw (7.45,1.0) node {$\cdots$};
\draw (7.4,1.3) node {$\overbrace{\phantom{xxxxx}}$};
\draw (7.4,1.8) node {$s-1$};
\draw (4.5,0.5) node {$\cdots$};
\end{tikzpicture}
\end{center}
\caption{Caterpillar $T(s,t)$}
\label{fig:st-graph}
\end{figure}

To see that Theorem~\ref{thm:lower-for-trees} is asymptotically optimal,
consider the caterpillars $T(s,t)$ shown in Figure~\ref{fig:st-graph}.

Clearly, $T(s,t)$ has $st$ vertices. Let $s\ge t+1$.  It is easy to
see that $\gamma_g(T(s,t)) = 2t-1$. Indeed, since $s-1\ge t$, \St can select a
leaf after each of the first $t-1$ moves of Dominator. Hence after \D selects
the $t$ vertices of high degree, the game is over.
By Theorem~\ref{thm:lower-for-trees},
$\gamma_g(T(s,t)) \ge \frac{2st}{s+4}-1$, which for a fixed $t$ and $n=st$
tends to $\frac{2n}{\Delta(T(s,t))+3} - 1 = \frac{2st}{s+4} -1 \sim 2t-1$ as $s\to \infty$.

\section{Pairs realizable by trees}
\label{sec:pairs}

In this section we are interested in which of the possible realizable pairs $(r,r)$,
$(r,r+1)$, and $(r,r-1)$ can be realized by trees. It was observed in~\cite{brklra-2010}
that $(k,k)$ is realizable by a tree, for $k\ge 1$. We now show that pairs
$(k,k+1)$ are also realizable by trees. On the other
hand, we prove that the pairs $(3,2)$ and $(4,3)$ cannot be realized by trees
and conjecture that no pair of the form $(k+1,k)$ is realizable by a tree.
(Clearly, no graph realizes the pair $(2,1)$.)

\begin{thm}
\label{thm:realize-plus1}
For any positive integer $k$, there exists a tree $T$ such that $\gamma_g(T)=k$ and $\gamma_g'(T)=k+1$.
\end{thm}

\proof
Stars confirm the result for $k=1$. For $k=2$ consider the
tree on five vertices obtained from $K_{1,3}$ by subdividing one edge.
In the rest of the proof assume that $k\ge 3$. We distinguish three cases based
on the congruence class of $k \pmod 3$.

\begin{figure}[ht!]
\begin{center}
\begin{tikzpicture}[scale=1.0,style=thick]
\def\vr{2pt} 
\draw (5.0,0.5) -- (0.0,2.0);
\draw (5.0,0.5) -- (3.0,2.0);
\draw (5.0,0.5) -- (10.0,2.0);
\draw (-1.0,2.0) -- (1.0,2.0);
\draw (2.0,2.0) -- (4.0,2.0);
\draw (9.0,2.0) -- (11.0,2.0);
\draw (5.0,0.5)  [fill=white] circle (\vr);
\draw (0.0,2.0)  [fill=white] circle (\vr);
\draw (3.0,2.0)  [fill=white] circle (\vr);
\draw (10.0,2.0)  [fill=white] circle (\vr);
\draw (-1.0,2.0)  [fill=white] circle (\vr);\draw (-0.5,2.0)  [fill=white] circle (\vr);
\draw (0.5,2.0)  [fill=white] circle (\vr);\draw (1.0,2.0)  [fill=white] circle (\vr);
\draw (2.0,2.0)  [fill=white] circle (\vr);\draw (2.5,2.0)  [fill=white] circle (\vr);
\draw (3.5,2.0)  [fill=white] circle (\vr);\draw (4.0,2.0)  [fill=white] circle (\vr);
\draw (9.0,2.0)  [fill=white] circle (\vr);\draw (9.5,2.0)  [fill=white] circle (\vr);
\draw (10.5,2.0)  [fill=white] circle (\vr);\draw (11.0,2.0)  [fill=white] circle (\vr);
\draw (5.0,0.2) node {$w$};
\draw (5.5,2.0) node {$\cdots$};
\draw (-0.5,2.3) node {$a_1$}; \draw (0.0,2.3) node {$b_1$}; \draw (0.5,2.3) node {$c_1$}; \draw (0.0,2.8) node {$X_1$};
\draw (2.5,2.3) node {$a_2$}; \draw (3.0,2.3) node {$b_2$}; \draw (3.5,2.3) node {$c_2$}; \draw (3.0,2.8) node {$X_2$};
\draw (9.5,2.3) node {$a_r$}; \draw (10.0,2.3) node {$b_r$}; \draw (10.5,2.3) node {$c_r$}; \draw (10.0,2.8) node {$X_r$};
\end{tikzpicture}
\end{center}
\caption{Tree $T_r$}
\label{fig:3realizations}
\end{figure}
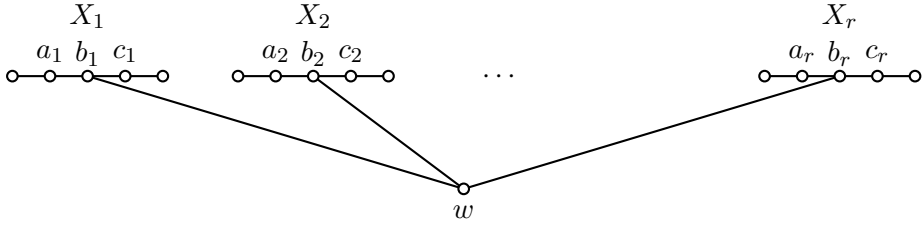

\medskip\noindent
{\bf Case 1}: $(3r,3r+1)$. \\
Let $r\ge 1$ and consider the tree $T_r$ of order $5r+1$ from Figure~\ref{fig:3realizations}.

\noindent
We will prove that $\gamma_g'(T_r) \ge 3r+1$ and $\gamma_g(T_r)\le 3r$.  These
two inequalities together with Theorem~\ref{thm:we-and-bill} show that $T_r$
realizes $(3r,3r+1)$.  We give a strategy for \St that will force at least
$3r+1$ vertices to be played in $T_r$.  Staller begins by playing $w$.
Her strategy is now to play in such a way that $b_t$ is played for every $t$ such
that $1 \le t \le r$.

She can accomplish this as follows. Dominator's first move will
be to play a vertex from some $X_i$.  By the Continuation Principle we see that
$d_1 \in \{a_i,b_i,c_i\}$.  If \D plays $a_i$ or $c_i$, then \St plays $b_i$.
On the other hand, if $d_1=b_i$, then \St plays $a_i$.
If \D now plays his second  move in $X_i$, then \St plays
$b_j$ for some $j$ different from $i$.  Otherwise, if \D plays his second move in $X_t$
for  $t\neq i$, then \St plays in $X_t$ using the same approach as she did in responding
to Dominator's first move in $X_i$.  By continuing this strategy \St can ensure that all
of the vertices in the set $\{b_1, \ldots,b_r\}$ are played in the course of the game.
This guarantees that at least $3r+1$ moves will be made, and thus $\gamma_g'(T_r) \ge 3r+1$.

Now consider Game 1 on $T_r$. \D begins by playing $a_1$. By symmetry and the
Continuation Principle, \St can choose essentially five different vertices for $s_1$.
These are $b_1$, $b_2$, $w$, the leaf $u$ adjacent to $c_1$, and the leaf $v$ adjacent to $a_2$.
For ease of explanation, let $P_3'$ denote
a partially dominated path of order 3 where one of the leaves is dominated, and let
$P_5'$  denote a partially dominated path of order 5 with the center vertex dominated. It is easy to see
that $\gamma_g(P_3')=1$, $\gamma_g'(P_3')=2$, and $\gamma_g(P_5')=3=\gamma_g'(P_5')$.

If $s_1=b_1$, then the residual graph after these
two moves is the disjoint union of two partially dominated trees, a path of order 2
with one of its vertices dominated and $T_{r-1}$ with $w$ dominated.  In this case it
follows by the Continuation Principle and induction that at most $3r$ moves will be made altogether.

If $s_1=w$, then the residual graph  is the disjoint union of $P_3'$ and $r-1$
copies of $P_5'$.  \D responds with $c_1$ in $P_3'$, and after that at most
$3(r-1)$ more vertices will be played.  If \St plays $u$ on her first move, then \D
responds with $w$.  The residual graph is now the disjoint union of $r-1$
copies of $P_5'$, and once again we see that at most $3r$ vertices will be played.  If $s_1=b_2$, then
\D plays $c_1$.  By this  move \D has limited the number of vertices played in
$X_1$ to 2 and can then play in such a way to ensure that no more than three vertices
are played from any $X_j$ with $j>1$.  The vertex $w$ might be played in the remainder
of the game, but we see again that a total of at most $3r$ moves will be
made in the game.

Finally, assume that $s_1=v$.  \D responds with $w$.
In this case the residual graph $F$ after these three moves is the disjoint
union of two copies of $P_3'$ and $r-2$ copies of $P_5'$.  Regardless of where \St
plays her second move, \D  plays in one of the copies of $P_3'$.  Since no additional move
on it will be played, it follows that on the corresponding $P_5$ only two moves are
played in the course of the game.  This ensures  that at most $3(r-1)$ vertices will be played in $F$, and again the total number of moves
in the game is no more than $3r$.  In all cases \D can limit the total number of
vertices played to $3r$, and hence $\gamma_g(T_r) \le 3r$.

As we noted at the beginning, it now follows that $T_r$ realizes $(3r,3r+1)$.

\medskip\noindent
{\bf Case 2}: $(3r+1,3r+2)$. \\
For $r\ge 1$ let $T_r'$ be the graph of order $5r+3$ obtained from $T_r$ (the
tree from Figure~\ref{fig:3realizations}) by
attaching a path of length 2 to $w$ with new vertices $y$ and $z$, where $z$
is a pendant vertex and $y$ is adjacent to $w$ and $z$.
Proceeding as we did above in Case 1, we show that $\gamma_g'(T_r') \ge 3r+2$ and $\gamma_g(T_r')\le 3r+1$.
Theorem~\ref{thm:we-and-bill} along with these two inequalities then imply that $T_r'$
realizes $(3r+1,3r+2)$.
In Game 2 \St first plays $w$, which leaves a residual graph that is the disjoint union
of $r$ copies of $P_5'$ and a path of order two with one of its vertices dominated.
Since \St can force at least three vertices to be played from each $P_5'$, it follows
that at least $3r+1$ more moves will be made on this residual graph, and hence
$\gamma_g'(T_r') \ge 1+ (3r+1)=3r+2$.

To begin Game 1 on $T_r'$, \D plays $a_1$.  Using symmetry and the Continuation Principle,
we conclude that \St has six different vertices to play as her first move.  That is, we
may assume $s_1\in\{b_1,u,w,z,b_2,v\}$, where $u$ and $v$ are the vertices of degree 1 as
described in Case 1.  If $s_1=b_1$, \D plays $a_2$; if $s_1\in\{u,v\}$, then \D responds with
$y$; if $s_1\in\{w,z,b_2\}$, then \D plays $c_1$.  With this second move by Dominator,
he can limit the total number of moves in Game 1 to at most $3r+1$.  The proof of this in
the six different cases is too detailed to include here, but it is similar to our analysis of
Game 1 in Case 1.  It now follows that $\gamma_g(T_r') \le 3r+1$, and thus $T_r'$ realizes
$(3r+1,3r+2)$.

\medskip\noindent
{\bf Case 3}: $(3r+2,3r+3)$. \\
In this case let $T_r''$ be the tree obtained from
the tree $T_r$ (of Figure~\ref{fig:3realizations}) by attaching two paths
of length 2 to $b_1$, say $P=b_1, p, q$  and $Q=b_1, m, n$.
We denote by $F$ the (partially dominated) subtree of $T_r''$ of order nine that is the
component of $T_r''-wb_1$ that contains $b_1$ and in which $b_1$ is dominated.
It can be shown that $F$ realizes $(5,5)$ and that $T_1''$ realizes $(5,6)$.
Because of the latter fact we assume that $r \ge 2$.  Let $S$ denote the component of $T_r''-wb_1$
that contains vertex $w$.

In Game 2 \St first plays $w$, leaving the residual graph $F \cup (r-1)P_5'$.
It follows that $\gamma_g'(T_r'') \ge 1+5+3(r-1)=3r+3$.

\D begins Game 1 on $T_r''$ by playing $a_2$.  After this opening move,
Dominator's goal is to play in such a way as either to prevent the vertex $w$
from being played in the course of the game or to limit the number of moves made on
some $P_5$ to 2. (By ``some $P_5$'' here we mean an induced path of order 5 that contains
both of $a_j$ and $c_j$ for some $j$ with $2 \le j \le r$.) Suppose that $s_1=w$.
Dominator then plays $c_2$, and the resulting residual graph is $F \cup (r-2)P_5'$.
In this case a total of at most $1+1+1+5+3(r-2)=3r+2$ moves will be made in the game.

Suppose instead that $s_1=b_2$.  In this case \D responds with $a_3$.   If
\St follows Dominator's moves by playing on the same $P_5$, then \D
continues to play $a_j$ from some $P_5$ that has not had one of its vertices played.
If, at some point in the game,  Staller plays $w$ before all of $c_2, c_3, \ldots, c_r$
are dominated, then \D can achieve his
goal by playing a second vertex on the (same) $P_5$  where he made his previous move to dominate
it in two moves.  Otherwise, \St will be the last player to play on $S$.  In this
case, \D plays the vertex $a_1$ thereby preventing the vertex $w$ from ever being played.
(It is easy to see that on any move by Staller in $F$ Dominator can follow in $F$
in such a way that the total of at most 5 moves will be played in $F$ during the game.)
Therefore, in all cases \D can ensure that at most $3r+2$ total moves are made in Game 1.

Again, we employ Theorem~\ref{thm:we-and-bill} to conclude that $T_r''$ realizes
$(3r+2,3r+3)$.
\qed

For the $(k,k-1)$ case we pose:

\begin{conj}
No pair of the form $(k,k-1)$ can be realized by a tree.
\end{conj}

In the rest of this section we prove the first two cases of the conjecture:

\begin{thm}
No tree realizes the pair $(3,2)$ or the pair $(4,3)$.
\end{thm}

\proof
Suppose that a tree $T$ realizes $(3,2)$. It is easy to see that
$\gamma_g'(T) = 2$ implies that $T$ is either a star $K_{1,n}$ for $n\ge 2$
or a $P_4$. In both cases $\gamma_g(T) \le 2$, so $(3,2)$ is not
realizable on trees.

Suppose $T$ is a tree that realizes $(4,3)$, and let $x$ be an optimal
first move for Dominator.  The residual graph $T'$ has at most 3
components, each of which is a partially dominated subtree of $T$.
Note that if one of these partially dominated components $F$ has
$\gamma_g'(F)=1$, then $F$ has exactly one undominated vertex.

Suppose first that $T'$ has three partially dominated components $T_1, T_2, T_3$
with $T_i$ rooted at the dominated vertex $v_i$.  If at least one of these
subtrees, say $T_1$, has more than one undominated vertex, then Staller
can force at least two moves in $T_1$.  Because the other two subtrees each
require at least one move, it follows that $\gamma_g(T) \ge 5$, a contradiction.
Hence, each of $T_1, T_2, T_3$ has exactly one undominated vertex, and $T$
is a tree formed by identifying a leaf from three copies of $P_3$ and
attaching some pendant vertices at the vertex of high degree.
However, this tree has Staller-start game domination number at least 4,
again contradicting our initial assumption.

Now suppose that $T'$ is the disjoint union of $T_1$ and $T_2$.  Note
that in this case $x$ cannot be a neighbor of a leaf in the original tree $T$.
Indeed, if $x$ is adjacent to a leaf $y$, then when Game 2 is played
on $T$, Staller can play first at $y$, which is easily shown to force at
least four moves.  Thus, $\deg(x)=2$.  If $\gamma_g'(T_1)=1=\gamma_g'(T_2)$,
then $T=P_5$ and $\gamma_g(T)=3$, a contradiction.

Note that the Staller-start game domination number of either of these two
partially dominated trees cannot exceed 2. We may thus assume without loss
of generality that $2=\gamma_g'(T_1) \ge \gamma_g'(T_2)$. Suppose that
$\gamma_g'(T_2)=2$. Staller can then play at vertex $x$ when Game 2 is played on $T$.
After Dominator's first move at least one of $T_1$ or $T_2$
is part of the residual graph, and Staller can then force at least two more
moves, again contradicting the assumption that $\gamma_g'(T)=3$.
Therefore, $T_2$ is the path of order 2 with one of its vertices dominated.

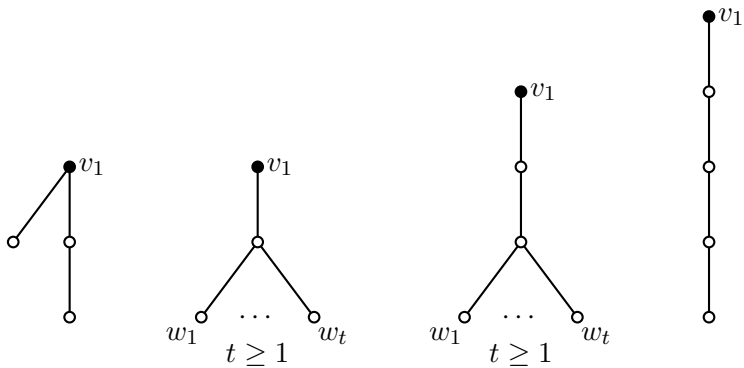
\begin{figure}[ht!]
\begin{center}
\begin{tikzpicture}[scale=1.0,style=thick]
\def\vr{2pt} 
\draw (1.0,0.0) -- (1.0,1.0) -- (1.0,2.0) -- (0.25,1.0);
\draw (2.75,0.0) -- (3.5,1.0) -- (3.5,2.0);
\draw (3.5,1.0) -- (4.25,0.0);
\draw (6.25,0.0) -- (7.0,1.0) -- (7.0,2.0) -- (7.0,3.0);
\draw (7.75,0.0) -- (7.0,1.0);
\draw (9.5,0.0) -- (9.5,4.0);
\draw (1.0,0.0)  [fill=white] circle (\vr);
\draw (1.0,1.0)  [fill=white] circle (\vr);
\draw (1.0,2.0)  [fill=black] circle (\vr);
\draw (0.25,1.0)  [fill=white] circle (\vr);
\draw (3.5,1.0)  [fill=white] circle (\vr);
\draw (3.5,2.0)  [fill=black] circle (\vr);
\draw (2.75,0.0)  [fill=white] circle (\vr);
\draw (4.25,0.0)  [fill=white] circle (\vr);
\draw (7.0,1.0)  [fill=white] circle (\vr);
\draw (7.0,2.0)  [fill=white] circle (\vr);
\draw (6.25,0.0)  [fill=white] circle (\vr);
\draw (7.75,0.0)  [fill=white] circle (\vr);
\draw (7.0,3.0)  [fill=black] circle (\vr);
\ fourth
\draw (9.5,0.0)  [fill=white] circle (\vr);
\draw (9.5,1.0)  [fill=white] circle (\vr);
\draw (9.5,2.0)  [fill=white] circle (\vr);
\draw (9.5,3.0)  [fill=white] circle (\vr);
\draw (9.5,4.0)  [fill=black] circle (\vr);

\draw (1.3,2.0) node {$v_1$};
\draw (3.8,2.0) node {$v_1$};
\draw (7.3,3.0) node {$v_1$};
\draw (9.8,4.0) node {$v_1$};
\draw (3.5,0.0) node {$\cdots$};
\draw (7.0,0.0) node {$\cdots$};
\draw (2.50,-0.25) node {$w_1$};
\draw (4.5,-0.25) node {$w_t$};
\draw (3.5,-0.5) node {$t\ge 1$};
\draw (6.00,-0.25) node {$w_1$};
\draw (8.0,-0.25) node {$w_t$};
\draw (7.0,-0.5) node {$t\ge 1$};

\end{tikzpicture}
\end{center}
\caption{Possible partially dominated trees}
\label{fig:4trees}
\end{figure}

If $T_1$ is a star with $v_1$ as its center or as one of its leaves, then
$\gamma(T)=2$ and hence $4=\gamma_g(T)\le 2\cdot 2 -1$, an obvious contradiction.
Therefore, $\gamma_g'(T_1)=2$, but $T_1$ is not a star.  A short analysis
shows that $T_1$ must be one of the partially dominated trees in Figure~\ref{fig:4trees}.
Each of these candidates for $T_1$ together with $T_2=P_2$ yields a tree
$T$ with either $\gamma_g(T)\not=4$ or $\gamma_g'(T)\not=3$, again contradicting
our assumption about $T$.  This implies that the residual graph $T'$ has
exactly one component.

Hence we are left with a tree $T$, a vertex $x$ that is an optimal
first move for Dominator, and the residual tree $T'$ which has one component.
Besides the neighbor $v_1$ in $T'$, the vertex $x$ is adjacent to some
leaves $y_1, \ldots, y_k$. We may assume that $k\ge 1$, because otherwise
(by the Continuation Principle) Dominator would rather select $v_1$ than
$x$ in his first move. Since $x$ is an optimal first move by \D in Game 1,
it follows that $\gamma_g'(T_1)=3$ in addition to $\gamma_g'(T)=3$.

Consider Game 2 played on the partially dominated tree $T_1$.
Let $w$ be an optimal first move by \St in this game, and let  $v$ in $T_1$
be an optimal response by Dominator.  At least one vertex, say  $u$, in $T_1$
is not dominated by $\{w,v\}$.  Note that $u \neq v_1$, since $T_1$ is a partially dominated tree
with $v_1$ dominated.  We can now show that $\gamma_g'(T) \ge 4$.  \St starts Game 2 on the original
tree $T$ by making the move  $s_1=w$.  Dominator either plays $d_1$ in $T_1$ or $d_1=x$.
If \D responds in $T_1$, then $y_1$ and at least one vertex in $T_1$ (other than $v_1$)
are not yet dominated.  Thus in this case at least four total moves are required in Game
2.  On the other hand, if $d_1=x$, then \St plays  $s_2=v$, and $u$ is not yet dominated.
Again Game 2 lasts at least a total of four moves.
This now implies that $\gamma_g'(T) \ge 4$, a contradiction.
\qed

\section{Game on spanning subgraphs}

We now turn our attention to relations between the game domination
number of a graph and its spanning subgraphs, in particular spanning
trees.

Note that since any graph is a spanning subgraph of the complete graph
of the same order, the ratio $\gamma_g(H)/\gamma_g(G)$ where $H$ is a
spanning subgraph of $G$ can be arbitrarily large.
On the other hand the following result shows that this ratio is bounded below by
one half.

\begin{prp}
If $G$ is a graph and $H$ is a spanning subgraph of $G$, then
$$\gamma_g(H) \ge \frac{\gamma_g(G) + 1}{2}\,.$$
In particular, if $T$ is a spanning tree of connected $G$, then
$\gamma_g(T) \ge (\gamma_g(G) + 1)/2$.
\end{prp}

\proof
Clearly, $\gamma(H) \ge \gamma(G)$. By Theorem~\ref{thm:connection},
$\gamma_g(H) \ge \gamma(H)$ and $\gamma_g(G) \le  2\gamma(G)-1$. Then
$\gamma_g(H) \ge \gamma(H) \ge \gamma(G) \ge (\gamma_g(G) + 1)/2$.
\qed

\begin{figure}[ht!]
\begin{center}
\begin{tikzpicture}[scale=0.8,style=thick]
\def\vr{2pt} 
\path (0,0) coordinate (Q);
\path(160:3cm) coordinate (P1);
\path(150:3.5cm) coordinate (P2);
\path(140:3.0cm) coordinate (P3);
\path(140:2.0cm) coordinate (P4);
\path(120:3cm) coordinate (R1);
\path(110:3.5cm) coordinate (R2);
\path(100:3.0cm) coordinate (R3);
\path(100:2.0cm) coordinate (R4);
\path(50:3cm) coordinate (S1);
\path(40:3.5cm) coordinate (S2);
\path(30:3.0cm) coordinate (S3);
\path(30:2.0cm) coordinate (S4);
\draw (Q) -- (P1) (P1) -- (P2) (P2) -- (P3) (Q) -- (P4) (P4) -- (P3) (P1) -- (P3);
\draw (Q) -- (R1) (R1) -- (R2) (R2) -- (R3) (Q) -- (R4) (R4) -- (R3) (R1) -- (R3);
\draw (Q) -- (S1) (S1) -- (S2) (S2) -- (S3) (Q) -- (S4) (S4) -- (S3) (S1) -- (S3);
\draw (Q)  [fill=white] circle (\vr);
\draw (P1)  [fill=white] circle (\vr); \draw (P2)  [fill=white] circle (\vr);
\draw (P3)  [fill=white] circle (\vr); \draw (P4)  [fill=white] circle (\vr);
\draw (R1)  [fill=white] circle (\vr); \draw (R2)  [fill=white] circle (\vr);
\draw (R3)  [fill=white] circle (\vr); \draw (R4)  [fill=white] circle (\vr);
\draw (S1)  [fill=white] circle (\vr); \draw (S2)  [fill=white] circle (\vr);
\draw (S3)  [fill=white] circle (\vr); \draw (S4)  [fill=white] circle (\vr);
\draw (0.2,2)  [fill=black] circle (1pt);
\draw (0.4,2)  [fill=black] circle (1pt);
\draw (0.6,1.95)  [fill=black] circle (1pt);
\path (8,0) coordinate (QQ);
\path(160:3cm) ++(8,0) coordinate (PP1);
\path(150:3.5cm) ++(8,0) coordinate (PP2);
\path(140:3.0cm) ++(8,0) coordinate (PP3);
\path(140:2.0cm) ++(8,0) coordinate (PP4);
\path(120:3cm) ++(8,0) coordinate (RR1);
\path(110:3.5cm) ++(8,0) coordinate (RR2);
\path(100:3.0cm) ++(8,0) coordinate (RR3);
\path(100:2.0cm) ++(8,0) coordinate (RR4);
\path(50:3cm) ++(8,0) coordinate (SS1);
\path(40:3.5cm) ++(8,0) coordinate (SS2);
\path(30:3.0cm) ++(8,0) coordinate (SS3);
\path(30:2.0cm) ++(8,0) coordinate (SS4);
\draw (QQ) -- (PP1) (PP1) -- (PP2) (PP2) -- (PP3) (QQ) -- (PP4)  (PP1) -- (PP3);
\draw (QQ) -- (RR1) (RR1) -- (RR2) (RR2) -- (RR3) (QQ) -- (RR4)  (RR1) -- (RR3);
\draw (QQ) -- (SS1) (SS1) -- (SS2) (SS2) -- (SS3) (QQ) -- (SS4)  (SS1) -- (SS3);
\draw (QQ)  [fill=white] circle (\vr);
\draw (PP1)  [fill=white] circle (\vr); \draw (PP2)  [fill=white] circle (\vr);
\draw (PP3)  [fill=white] circle (\vr); \draw (PP4)  [fill=white] circle (\vr);
\draw (RR1)  [fill=white] circle (\vr); \draw (RR2)  [fill=white] circle (\vr);
\draw (RR3)  [fill=white] circle (\vr); \draw (RR4)  [fill=white] circle (\vr);
\draw (SS1)  [fill=white] circle (\vr); \draw (SS2)  [fill=white] circle (\vr);
\draw (SS3)  [fill=white] circle (\vr); \draw (SS4)  [fill=white] circle (\vr);
\draw (8.2,2)  [fill=black] circle (1pt);
\draw (8.4,2)  [fill=black] circle (1pt);
\draw (8.6,1.95)  [fill=black] circle (1pt);
\draw (2.3,0.5) node {$G_t$};
\draw (10.3,0.5) node {$H_t$};
\draw (-2.0,1.2) node {$1$};
\draw (-0.8,2.0) node {$2$};
\draw (1.7,1.5) node {$t$};
\draw (0.3,-0.3) node {$x$};
\draw (8.3,-0.3) node {$x$};
\end{tikzpicture}
\end{center}
\caption{Graph $G_t$ and its spanning subgraph $H_t$}
\label{fig:spanning-graphs}
\end{figure}
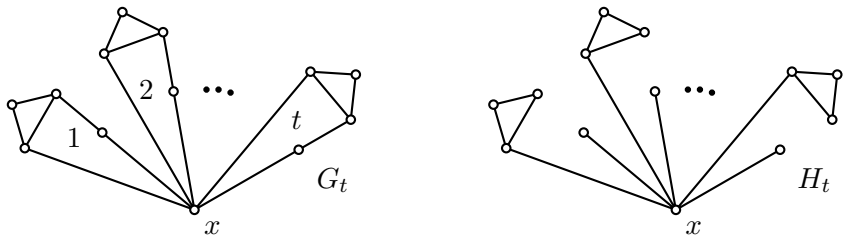

To see that a spanning subgraph can indeed have game domination number much smaller
than its supergraph, consider the graph $G_t$ consisting of $t$ blocks isomorphic
to $\overline{P_5}$ (the graph obtained from $C_5$ by adding an edge) and its spanning subgraph
$H_t$, see Figure~\ref{fig:spanning-graphs}.
Let $x$ be the vertex where the houses of $G_t$ are amalgamated. Note that
\D needs at least two moves to dominate each of the blocks of $G_t$. Hence
his strategy is to play $x$ and then finish dominating one block
on each move. On the other hand, if not all blocks are already dominated, \St
can play the vertex of degree 2 adjacent to $x$ of such a block $B$ in order
to force one more move on $B$. So in half of the blocks two vertices will be
played (not counting the move on $x$), which in turn implies that $\gamma_g(G_t)$
is about $3t/2$. On the other hand, playing Game 1 on $H_t$, the optimal first move
for \D is $x$. After that \St and \D will in turn dominate each of the
triangles, so $\gamma_g(H_t) = t+1$.

The example of Figure~\ref{fig:spanning-graphs} might lead to a suspicion that
2-connected spanning subgraphs cannot have smaller game domination number than their
2-connected supergraphs. However:

\begin{thm}
For $m\ge 3$, there exists a 3-connected graph $G_m$ having a 2-connected
spanning subgraph $H_m$ such that $\gamma_g(G_m) \ge 2m-2$ and $\gamma_g(H_m) = m$.
\end{thm}

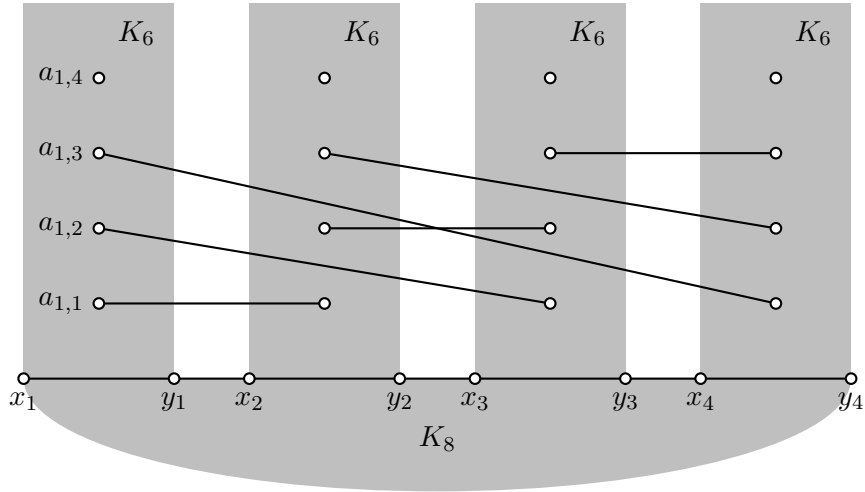
\begin{figure}[ht!]
\begin{center}
\begin{tikzpicture}[scale=1.0,style=thick]
\def\vr{2pt} 

\fill[color=lightgray] (0.0,5.0) rectangle (2.0,0.0);
\fill[color=lightgray] (3.0,5.0) rectangle (5.0,0.0);
\fill[color=lightgray] (6.0,5.0) rectangle (8.0,0.0);
\fill[color=lightgray] (9.0,5.0) rectangle (11.0,0.0);

\begin{scope}
\clip (0.0,0.0) rectangle (11.0,-1.5);
\clip (5.5,0.0) ellipse (5.5 and 1.5);
\fill[color=lightgray] (0.0,0.0) rectangle (11.0,-3.0);
\end{scope}

\path (0,0) coordinate (v1); \path (2,0) coordinate (v2);
\path (3,0) coordinate (v3); \path (5,0) coordinate (v4);
\path (6,0) coordinate (v5); \path (8,0) coordinate (v6);
\path (9,0) coordinate (v7); \path (11,0) coordinate (v8);
\path (1,1) coordinate (x1); \path (1,2) coordinate (x2);
\path (1,3) coordinate (x3); \path (1,4) coordinate (x4);
\draw (v1) -- (v8);
\draw (x1) -- (4.0,1.0); \draw (x2) -- (7.0,1.0);
\draw (x3) -- (10.0,1.0);
\draw (4.0,2.0) -- (7.0,2.0); \draw (4.0,3.0) -- (10.0,2.0);
\draw (7.0,3.0) -- (10.0,3.0);
\foreach \i in {1,...,8}
{  \draw (v\i)  [fill=white] circle (\vr); }
\foreach \i in {1,...,4}
{  \draw (x\i)  [fill=white] circle (\vr); }
\foreach \i in {1,...,4}
{  \draw (x\i) ++(3.0,0.0)  [fill=white] circle (\vr); }
\foreach \i in {1,...,4}
{  \draw (x\i) ++(6.0,0.0)  [fill=white] circle (\vr); }
\foreach \i in {1,...,4}
{  \draw (x\i) ++(9.0,0.0)  [fill=white] circle (\vr); }
\draw (5.5,-0.8) node {$K_8$};
\draw (1.5,4.6) node {$K_6$};
\draw (4.5,4.6) node {$K_6$};
\draw (7.5,4.6) node {$K_6$};
\draw (10.5,4.6) node {$K_6$};
\draw (0.0,-0.3) node {$x_1$};
\draw (2.0,-0.3) node {$y_1$};
\draw (3.0,-0.3) node {$x_2$};
\draw (5.0,-0.3) node {$y_2$};
\draw (6.0,-0.3) node {$x_3$};
\draw (8.0,-0.3) node {$y_3$};
\draw (9.0,-0.3) node {$x_4$};
\draw (11.0,-0.3) node {$y_4$};
\draw (0.5,1.0) node {$a_{1,1}$};
\draw (0.5,2.0) node {$a_{1,2}$};
\draw (0.5,3.0) node {$a_{1,3}$};
\draw (0.5,4.0) node {$a_{1,4}$};
\end{tikzpicture}
\end{center}
\caption{Graph $G_4$}
\label{fig:3-connected}
\end{figure}

\proof
We form a graph $G_m$ of order $m(m+2)$ as follows. Let
$X_i= \{a_{i,1}, \ldots, a_{i,m} \} \cup \{x_i, y_i \}$ for
$1 \le i\le m$, and then set
\[ V(G_m)= \bigcup_{i=1}^m X_i\,.\]
The edges are the following.  Let
$\{x_1,y_1,\ldots, x_m,y_m\}$ induce a complete graph of
order $2m$.  For $1 \le p \le m$, let $X_i$ induce a complete graph of order $m+2$.
In addition, for $1\le i\le m-1$ and $i\le j\le m-1$, add the edge $a_{i,j}a_{j+1,i}$.
See Figure~\ref{fig:3-connected} for $G_4$.

If $d_1 = x_1$, then \St plays in $X_1$, say at $a_{1,1}$.
In each of the subsequent rounds, the Continuation Principle implies that
\D must play in some $X_i$ that has not been played in before and on a vertex
of $X_i$ that has an undominated neighbor outside $X_i$. It will always be possible for \St
to follow \D and also play in $X_i$ in each of her first $m-2$ moves.
Hence by this time, $2m-4$ moves are made. At this stage, there are two
undominated vertices in different $X_i$'s with no common neighbor. Hence
two more moves are needed to end the game, which thus ends in no less than $2m-2$ moves.

Assume next that $d_1 = a_{1,1}$.  Now \St plays $x_1$ and we are
in the first case. Note that $d_1=a_{1,m}$ need not be considered due to the
Continuation Principle, so $d_1\in \{x_1,a_{1,1}\}$ covers all cases by symmetry.
Hence $\gamma_g(G_m) \ge 2m-2$.

The spanning subgraph $H_m$ of $G_m$ is obtained by removing all the edges $a_{i,j}a_{j+1,i}$.
By the Continuation Principle, we may without loss of generality assume that
$d_1 = x_1$ when Game 1 is played on $H_m$. In this case, each successive move of either player
completely dominates the $X_i$ in which it is played. Hence $\gamma_g(H_m) = m$.
\qed

If $\gamma_g(G)$ attains one of the two possible extremal values,
$\gamma(G)$ or $2\gamma(G)-1$, then we can say more.

\begin{prp}
(i) If $G$ is a graph with $\gamma_g(G) = \gamma(G)$ and  $H$ is a spanning
subgraph of $G$, then $\gamma_g(H) \ge \gamma_g(G)$.

(ii) If $G$ is a graph with $\gamma_g(G) = 2\gamma(G)-1$ and  $H$ is a spanning
subgraph of $G$ with $\gamma(H) = \gamma(G)$, then $\gamma_g(H) \le \gamma_g(G)$.
\end{prp}

\proof
Assertion (i) follows because $\gamma_g(H) \ge \gamma(H) \ge \gamma(G) = \gamma_g(G)$;
while (ii) follows since $\gamma_g(H) \le 2\gamma(H)-1 = 2\gamma(G)-1 = \gamma_g(G)$.
\qed

Since every graph $G$ has a spanning forest $F$ such that $\gamma(G) = \gamma(F)$,
see~\cite[Exercise 10.14]{imklra-2008}, we also infer that if $G$ is a graph
with $\gamma_g(G) = 2\gamma(G)-1$, then $G$ contains a spanning forest $F$
(spanning tree if $G$ is connected) such that $\gamma_g(F) \le \gamma_g(G)$.

In the rest of this section we focus on spanning trees.
First we show that {\bf all} spanning trees of a graph $G$ may have
game domination number much larger than $G$.

\begin{thm}
For $m=2r$, $r\ge 1$, there exists a graph $G_m$ such that
$$\gamma_g(T) - \gamma_g(G_m) \ge r-1$$
holds for any spanning tree $T$ of $G_m$.
\end{thm}

\begin{figure}[ht!]
\begin{center}
\begin{tikzpicture}[scale=0.8,style=thick]
\def\vr{2pt} 
\draw (2.0,4.5)--(-2.0,2.0);
\draw (2.0,4.5)--(0.0,2.0);
\draw (2.0,4.5)--(5.0,2.0);
\draw (2.0,4.5)  [fill=white] circle (\vr);
\draw (-2.0,2.0)  [fill=white] circle (\vr);
\draw (0.0,2.0)  [fill=white] circle (\vr);
\draw (5.0,2.0)  [fill=white] circle (\vr);
\draw (2.0,4.8) node {$x$};
\draw (-2.2,2.4) node {$v_1$};
\draw (-0.2,2.4) node {$v_2$};
\draw (5.2,2.4) node {$v_m$};
\draw (-2.0,1.0) node {$K_n$};
\draw (0.0,1.0) node {$K_n$};
\draw (5.0,1.0) node {$K_n$};
\draw (2.5,3.0) node {$\cdots$};
\draw (-2.0,0.8) ellipse (20pt and 40pt);
\draw (0.0,0.8) ellipse (20pt and 40pt);
\draw (5.0,0.8) ellipse (20pt and 40pt);
\end{tikzpicture}
\end{center}
\caption{Graph $G_m$}
\label{fig:star-complete}
\end{figure}
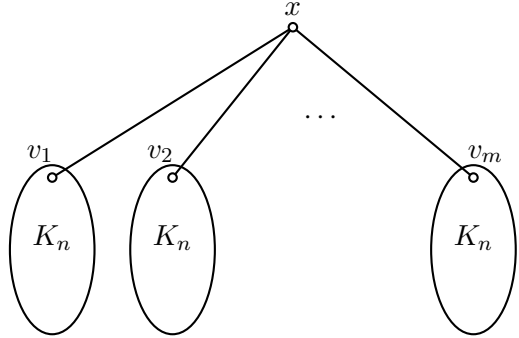

\proof
Let $n\ge 3$ and let $S$ be the star with center $x$ and
leaves $v_1, \ldots, v_m$. Let $G_m$ be the graph (of order $nm+1$) constructed
by identifying a vertex of a complete graph of order $n$ with $v_i$, for $1\le i\le m$;
see Figure~\ref{fig:star-complete}.

We first note that $\gamma_g(G_m) = m+1$.
Let $T$ be any spanning tree of $G_m$. $T$ has at least one leaf $\ell_i$ in the subtree
$T_i$ of $T$ rooted at $v_i$ when the edge $xv_i$ is removed from $T$ (choose $\ell_i \not= v_i$).
When Game 1 is played on $T$, \St can choose at least half of these leaves ($\ell_1,
\ldots, \ell_m$) or let \D choose them. Thus in at least half of $T_1, \ldots, T_m$,
two vertices will be chosen. Therefore $\gamma_g(T) \ge m + m/2 = 3m/2$, so
$$\gamma_g(T) - \gamma_g(G_m) \ge \frac{3}{2}m - m - 1 = r-1\,.$$
\qed

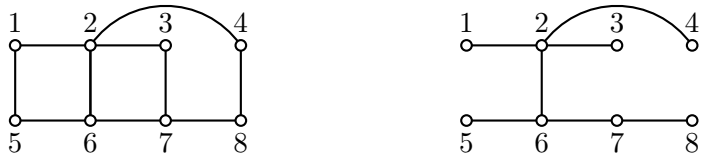
\begin{figure}[ht!]
\begin{center}
\begin{tikzpicture}[scale=1.0,style=thick]
\def\vr{2pt} 
\draw (0.0,0.0) rectangle (1.0,1.0);
\draw (1.0,0.0) rectangle (2.0,1.0);
\draw (2.0,0.0)--(3.0,0.0);
\draw (3.0,0.0)--(3.0,1.0);
\draw (1.0,1.0) .. controls (1.5,1.7) and (2.5,1.7) .. (3.0,1.0);
\draw (7.0,1.0) .. controls (7.5,1.7) and (8.5,1.7) .. (9.0,1.0);
\draw (6.0,0.0)--(9.0,0.0);
\draw (6.0,1.0)--(8.0,1.0);
\draw (7.0,0.0)--(7.0,1.0);
\draw (0.0,0.0)  [fill=white] circle (\vr);
\draw (1.0,0.0)  [fill=white] circle (\vr);
\draw (2.0,0.0)  [fill=white] circle (\vr);
\draw (3.0,0.0)  [fill=white] circle (\vr);
\draw (0.0,1.0)  [fill=white] circle (\vr);
\draw (1.0,1.0)  [fill=white] circle (\vr);
\draw (2.0,1.0)  [fill=white] circle (\vr);
\draw (3.0,1.0)  [fill=white] circle (\vr);
\draw (6.0,0.0)  [fill=white] circle (\vr);
\draw (7.0,0.0)  [fill=white] circle (\vr);
\draw (8.0,0.0)  [fill=white] circle (\vr);
\draw (9.0,0.0)  [fill=white] circle (\vr);
\draw (6.0,1.0)  [fill=white] circle (\vr);
\draw (7.0,1.0)  [fill=white] circle (\vr);
\draw (8.0,1.0)  [fill=white] circle (\vr);
\draw (9.0,1.0)  [fill=white] circle (\vr);
\draw (0.0,1.3) node {$1$};
\draw (1.0,1.3) node {$2$};
\draw (2.0,1.3) node {$3$};
\draw (3.0,1.3) node {$4$};
\draw (0.0,-0.3) node {$5$};
\draw (1.0,-0.3) node {$6$};
\draw (2.0,-0.3) node {$7$};
\draw (3.0,-0.3) node {$8$};
\draw (6.0,1.3) node {$1$};
\draw (7.0,1.3) node {$2$};
\draw (8.0,1.3) node {$3$};
\draw (9.0,1.3) node {$4$};
\draw (6.0,-0.3) node {$5$};
\draw (7.0,-0.3) node {$6$};
\draw (8.0,-0.3) node {$7$};
\draw (9.0,-0.3) node {$8$};
\end{tikzpicture}
\end{center}
\caption{Graph $G$ and its spanning tree $T$}
\label{fig:bad-st}
\end{figure}

Recall that the domination number of a spanning subgraph of a graph $G$ cannot be smaller than that of $G$.
In contrast we now give an example of a graph $G$ with a spanning tree $T$ such that $\gamma_g(T)<\gamma_g(G)$.
Consider the graph $G$ and the spanning tree $T$ from Figure~\ref{fig:bad-st}.

\begin{figure}[ht!]
\begin{center}
\begin{tikzpicture}[scale=1.3,style=thick]
\def\vr{2pt} 
\draw (0.5,0.0) -- (1.0,0.0);
\draw (1.0,0.0) -- (1.5,0.0);
\draw (3.0,0.0) -- (3.5,0.0);
\draw (5.0,0.0) -- (5.5,0.0);
\draw (7.0,0.0) -- (7.5,0.0);
\draw (0.5,0.0)  [fill=black] circle (\vr);
\draw (1.0,0.0)  [fill=white] circle (\vr);
\draw (1.5,0.0)  [fill=white] circle (\vr);
\draw (3.0,0.0)  [fill=black] circle (\vr);
\draw (3.5,0.0)  [fill=white] circle (\vr);
\draw (5.0,0.0)  [fill=white] circle (\vr);
\draw (5.5,0.0)  [fill=black] circle (\vr);
\draw (7.0,0.0)  [fill=white] circle (\vr);
\draw (7.5,0.0)  [fill=black] circle (\vr);

\draw (0.0,0.0) node {$(5,$};
\draw (1.7,0.0) node {$);$};
\draw (0.5,-0.3) node {$6$};
\draw (1.0,-0.3) node {$7$};
\draw (1.5,-0.3) node {$8$};
\draw (2.5,0.0) node {$(6,$};
\draw (3.7,0.0) node {$);$};
\draw (3.0,-0.3) node {$7$};
\draw (3.5,-0.3) node {$8$};
\draw (4.5,0.0) node {$(7,$};
\draw (5.7,0.0) node {$);$};
\draw (5.0,-0.3) node {$5$};
\draw (5.5,-0.3) node {$6$};
\draw (6.5,0.0) node {$(8,$};
\draw (7.7,0.0) node {$);$};
\draw (7.0,-0.3) node {$5$};
\draw (7.5,-0.3) node {$6$};

\end{tikzpicture}
\end{center}
\caption{Staller's possible moves}
\label{fig:on-st}
\end{figure}
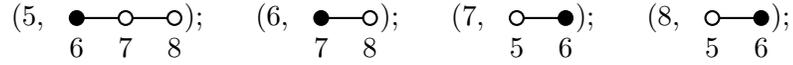

For each of the following pairs $(x,y)$ of vertices  from $G$, if \D plays $x$ then \St can play
$y$ and then the game domination number of the resulting residual graph $G'$
will be 2: $(1,6); (2,3); (3,2); (4,8); (8,4); (7,3); (6,1); (5,1)$. Therefore,
$\gamma_g(G)\ge 4$. Consider now the spanning tree $T$, and let \D play 2 on $T$. For each
of the following vertices $a$, the residual graph $T'$ when \St plays $a$ is listed
in Figure~\ref{fig:on-st}. For instance, the left case is when \St plays $5$; in this
case the residual graph is induced by vertices $6, 7, 8$ and the vertex $6$ of the residual
graph is already dominated, as indicated by the filled vertex.

In each case we find that the residual graph has game domination
number 1, so
$$\gamma_g(T) \le 3 < \gamma_g(G)\,.$$

This rather surprising fact demonstrates the intrinsic difficulty and
unusual behavior of the game domination number. Even more can be shown:

\begin{thm}
For any positive integer $\ell$, there exists a graph $G$, having a spanning
tree $T$ such that $\gamma_g(G)-\gamma_g(T)\ge \ell$.
\end{thm}

\begin{figure}[ht!]
\begin{center}
\begin{tikzpicture}[scale=0.8,style=thick]
\def\vr{2pt} 

\path (0.0,0.0) coordinate (y); \path (0.0,5.0) coordinate (x);
\path (-1,3.2) coordinate (x1); \path (-1.5,3.4) coordinate (x2);
\path (-2.0,3.6) coordinate (x3); \path (-2.5,3.8) coordinate (x4);
\path (-3.0,4.0) coordinate (x5); \path (-3.5,4.2) coordinate (x6);
\path (-4.0,4.4) coordinate (x7); \path (-4.5,4.6) coordinate (x8);
\path (-5.0,4.8) coordinate (x9); \path (-5.5,5.0) coordinate (x10);
\path (1,3.2) coordinate (w1); \path (1.5,3.4) coordinate (w2);
\path (2.0,3.6) coordinate (w3); \path (2.5,3.8) coordinate (w4);
\path (3.0,4.0) coordinate (w5); \path (3.5,4.2) coordinate (w6);
\path (4.0,4.4) coordinate (w7); \path (4.5,4.6) coordinate (w8);
\path (5.0,4.8) coordinate (w9); \path (5.5,5.0) coordinate (w10);
\path (-1,1.8) coordinate (y1); \path (-1.5,1.6) coordinate (y2);
\path (-2.0,1.4) coordinate (y3); \path (-2.5,1.2) coordinate (y4);
\path (-3.0,1.0) coordinate (y5); \path (-3.5,0.8) coordinate (y6);
\path (-4.0,0.6) coordinate (y7); \path (-4.5,0.4) coordinate (y8);
\path (-5.0,0.2) coordinate (y9); \path (-5.5,0.0) coordinate (y10);
\path (1,1.8) coordinate (z1); \path (1.5,1.6) coordinate (z2);
\path (2.0,1.4) coordinate (z3); \path (2.5,1.2) coordinate (z4);
\path (3.0,1.0) coordinate (z5); \path (3.5,0.8) coordinate (z6);
\path (4.0,0.6) coordinate (z7); \path (4.5,0.4) coordinate (z8);
\path (5.0,0.2) coordinate (z9); \path (5.5,0.0) coordinate (z10);
\draw (x) -- (y);
\foreach \i in {1,...,10}
{ \draw (x) -- (x\i);
  \draw (x) -- (w\i);
  }
\draw (y) -- (y5);
\draw (y) -- (y10);
\draw (y) -- (z5);
\draw (y) -- (z10);
\draw (y1) -- (y5);
\draw (y6) -- (y10);
\draw (z1) -- (z5);
\draw (z6) -- (z10);
\foreach \i in {1,...,10}
{ \draw (x\i) -- (y\i);
  \draw (w\i) -- (z\i);
  }

\draw (x)  [fill=white] circle (\vr);
\draw (y)  [fill=white] circle (\vr);
\foreach \i in {1,...,10}
{ \draw (x\i)  [fill=white] circle (\vr);
  \draw (w\i)  [fill=white] circle (\vr);
  \draw (y\i)  [fill=white] circle (\vr);
  \draw (z\i)  [fill=white] circle (\vr);
  }
\draw (0.0,-0.3) node {$y$};
\draw (0.0,5.3) node {$x$};
\end{tikzpicture}
\end{center}
\caption{Graph $G_4$}
\label{fig:cuden}
\end{figure}
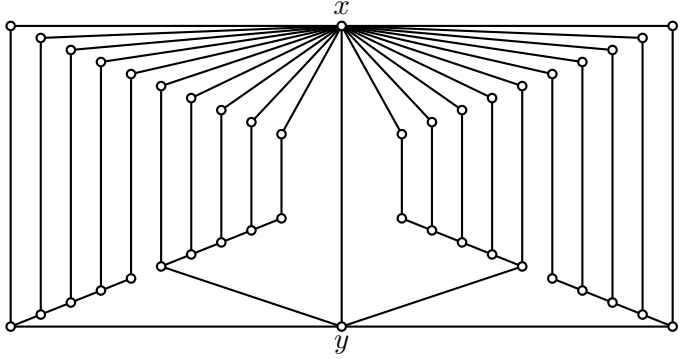

\proof
We introduce the family of graphs $G_k$ and their spanning trees $T_k$
in the following way.
Let $k$ be a positive integer, and for each $i$ between 1 and $k$,
$x_{i}^1,x_{i}^2,x_{i}^3,x_{i}^4,x_{i}^5$ are non-adjacent vertices in $T_k$,
and $Q_i:y_{i}^1y_{i}^2y_{i}^3y_{i}^4y_{i}^5$ is a path isomorphic to $P_5$ in $T_k$.
Finally $x$ and $y$ are two vertices, such that $x$ is adjacent to
$x_{i}^1,x_{i}^2,x_{i}^3,x_{i}^4$ and $x_{i}^5$ for all $i\in \{1,\ldots,k\}$,
while $y$ is adjacent to $y_{i}^1$ for all $i\in \{1,\ldots,k\}$,
and $x$ and $y$ are also adjacent. The resulting graph $T_k$ is a tree
on $10k+2$ vertices. We obtain $G_k$ by adding edges between $x_{i}^j$
and $y_i^j$ for $1\le i\le k$, $1\le j\le 5$.
See Figure~\ref{fig:cuden} for $G_4$, from which $T_4$ is obtained by
removing all vertical edges except $xy$.

To complete the proof it suffices to show that for any integer $k\geq 1$,
$$\gamma_g(G_k)\geq \frac{5}{2}k-1\quad {\rm and}\quad \gamma_g(T_k)\leq {2}k+3\,.$$
Let us first verify the second inequality, concerning trees $T_k$.
To prove it we need to show that \D has a strategy by which at
most $2k+3$ moves will be played during Game 1. His strategy is as
follows. In his first two moves, he ensures that $x$ and $y$ are
chosen. He plays $x$ in his first move, and $y$ in his second
move, unless already \St played $y$ (we will consider this case
later). Now, $s_1\ne y$ implies that $s_1$ is in some $Q_i$;
without loss of generality let this be $Q_1$. Hence in Dominator's third move,
since $y_i^1$ is dominated for each $i$, he can dominate all vertices of $Q_1$.
One by one, \St will have to pick a new $Q_i$ to play in, which \D will entirely
dominate in his next move. Altogether, in each $Q_i$ (with a possible
exception of one $Q_i$, where \St can force three vertices to be played),
there will be only two vertices played, which yields $2k+3$ as the total number
of moves in this game. On the other hand, if $s_1=y$, then $d_2=y_1^3$ ensures that in $Q_1$
only two vertices will be played. In addition, by a similar strategy as above \D can force
that only two moves will be played in each of $Q_i$s. Hence only $2k+2$ moves will be played.

To prove the first inequality we need to show
that \St has a strategy to enforce at least $\frac{5}{2}k-1$ moves played
during Game 1 in $G_k$. Her strategy in each of the first $k$ moves of the game is
to play an $x_i^4$ such that no vertex from $Q_i\cup \{ x_{i}^1, x_{i}^2, x_{i}^3, x_{i}^5\}$
has yet been played.
Using this strategy she ensures that at least two more moves will be needed to dominate
each of these $\lfloor \frac{k}{2}\rfloor$ $Q_i$s (since at least $y_i^2,y_i^3$ and $y_i^5$
are left undominated). The remaining paths $Q_i$
require at least two moves each as well. Hence altogether, there will be at least
$2k+\lfloor \frac{k}{2}\rfloor$ moves played during Game 1, which implies
$\gamma_g(G_k)\geq \frac{5}{2}k-1$.
\qed

\section*{Acknowledgements}
The authors extend their appreciation to anonymous referees for their
careful reading of our manuscript and for the resulting suggestions for improvement.
In particular, several of the proofs in the paper were clarified with their help.
We are also grateful to Editor-in-Chief Douglas B. West for investing considerable
time and effort in pointing out grammatical (and some other) corrections.


\end{document}